\documentclass[11pt]{article}
\usepackage[utf8]{inputenc}
\usepackage{mathtools}
\usepackage{amsfonts}
\usepackage{amssymb}
\usepackage{amsthm}
\usepackage{amsmath, amscd, amsfonts, amssymb}
\usepackage{enumitem}
\usepackage{authblk}
\newtheorem{theorem}{Theorem}
\newtheorem{corollary}[theorem]{Corollary}
\newtheorem{lemma}[theorem]{Lemma}
\newtheorem{question}[theorem]{Question}
\def\lc{\left\lceil}   
\def\rc{\right\rceil}
\def\lf{\left\lfloor}   
\def\rf{\right\rfloor}
\makeatletter
\def\blfootnote{\gdef\@thefnmark{}\@footnotetext}
\makeatother
\makeatletter
\@namedef{subjclassname@2020}{\textup{2020} Mathematics Subject Classification}
\makeatother

\title{Factoring complete graphs and hypergraphs into factors with few maximal cliques}

\author[1]{Paul Erdős}
\author[2]{David P. Galvin\thanks{david.galvin@gmail.com}}
\author[3]{Fred Galvin\thanks{bof@sunflower.com}}
\author[4]{Michael M. Krieger\thanks{mkrieger@alumni.caltech.edu}}
\affil[1]{Hungarian Academy of Sciences, Budapest, Hungary}
\affil[2]{Somerville, MA 02143}
\affil[3]{Dept. of Mathematics, University of Kansas, Lawrence, KS 66045}
\affil[4]{Krieger Law Office, 11209 National Blvd., \#417, Los Angeles, CA 90064}

\date{%
    \today
}

\begin{document}
\maketitle

\begin{abstract}
For integers $r,t\geq2$ and $n\geq1$ let $f_r(t,n)$ be the minimum, over all factorizations of the complete $r$-uniform hypergraph of order $n$ into $t$ factors $H_1,\dots,H_t$, of $\sum_{i=1}^tc(H_i)$ where $c(H_i)$ is the number of maximal cliques in $H_i$. It is known that $f_2(2,n)=n+1$; in fact, if $G$ is a graph of order $n$, then $c(G)+c(\overline G)\geq n+1$ with equality iff $\omega(G)+\alpha(G)=n+1$ where $\omega$ is the clique number and $\alpha$ the independence number. In this paper we investigate $f_r(t,n)$ when $r>2$ or $t>2$. We also characterize graphs $G$ of order $n$ with $c(G)+c(\overline G)=n+2$.
\end{abstract}

\blfootnote{Some of this research was conducted at the University of California, Los Angeles, and the University of Kansas, and was supported in part by grants from the National Science Foundation.}

\blfootnote{\textup{2020} \emph{Mathematics Subject Classification}:
Primary: 05C69; Secondary: 05C65}


\section{Introduction}

Let $c(G)$ denote the number of maximal cliques and $\overline c(G)=c(\overline G)$ the number of maximal anticliques in a graph $G$. In an earlier paper \cite{GK71} two of the present authors showed that, if $G$ is a graph of order $n$, then $c(G)+\overline c(G)\geq n+1$ with equality just in case $\omega(G)+\alpha(G)=n+1$ where $\omega$ is the clique number and $\alpha$ the independence number. (Proofs of these facts will be given here; see Theorem \ref{T28} and Theorem \ref{T34}(a).) In this paper we generalize the problem in two ways: by considering factorizations of the complete graph $K_n$ into $t$ factors instead of just the two factors $G$ and $\overline G$, and by considering uniform hypergraphs instead of graphs. We also characterize the graphs $G$ of order $n$ with $c(G)+\overline c(G)=n+2$.

An \emph{$r$-uniform hypergraph} (in this paper $r\ge2$) is a structure $H=(V,E)$ comprising a nonempty finite set $V=V(H)$ of \emph{vertices} and a set $E=E(H)\subseteq\binom Vr$ of \emph{edges}; $n(H)=|V|$ is the \emph{order} and $e(H)=|E|$ is the \emph{size} of $H$. (A $2$-uniform hypergraph is a \emph{graph}; a $3$-uniform hypergraph is a \emph{triple system}.) A set $X\subseteq V$ is a \emph{clique} if $\binom Xr\subseteq E$, an \emph{anticlique} if $\binom Xr\cap E=\varnothing$. The \emph{complement} of $H$ is the hypergraph $\overline H=(V,\binom Vr\setminus E)$. We write $c(H)$ for the number of maximal cliques and $\overline c(H)=c(\overline H)$ for the number of maximal anticliques in $H$. We write $D(v)$ for the number of maximal cliques and $\overline D(v)$ for the number of maximal anticliques containing a vertex $v$. We define $d(H)=\min\{D(v):v\in V\}$ and $\overline d(H)=d(\overline H)=\min\{\overline D(v):v\in V\}$. The \emph{complete} $r$-uniform hypergraph of order $n$ is $K^r_n=(V,\binom Vr)$ where $|V|=n$. Note that $c(K^r_n)=1$ and $\overline c(K^r_n)=\max\{\binom n{r-1},1\}$. An edge $\{u,v\}$ of a graph may be written as $uv$ if no confusion will result. The \emph{neighborhood} of a vertex $v$ in a graph $G=(V,E)$ is the set $N(v)=N_G(v)=\{u\in V:uv\in E\}$, and we write $\overline N(v)$ for $N_{\overline G}(v)=\{u\in V:uv\notin E,u\ne v\}$. The \emph{clique number} $\omega(G)$ of a graph $G$ is the maximum number of vertices in a clique of $G$; the \emph{independence number} $\alpha(G)=\omega(\overline G)$ is the maximum number of vertices in an anticlique of $G$. See West \cite{W01} for graph-theoretic notation and terminology not defined here.

For integers $r,t\geq2$ and $n\geq1$ we define $f_r(t,n)$ as the minimum of $\sum_{i=1}^tc(H_i)$ over all factorizations of $K^r_n$ into $t$ edge-disjoint spanning subgraphs $H_1,\dots,H_t$. In this notation, part of the aforementioned result of \cite{GK71} may be expressed as follows. (A slightly stronger result will be proved here as Theorem \ref{T28}.)

\begin{theorem}
\label{T1}
$f_2(2, n) = n + 1$.
\end{theorem}

Sometimes we find it convenient to use the language of colorings: $f_r(t,n)$ is the least possible number of maximal monochromatic cliques in an edge coloring of $K^r_n$ with $t$ colors, it being understood that a set which is a maximal monochromatic clique for more than one color is counted once for each color; e.g., $f_2(3,2)=1+2+2=5$. Most of this paper is a study of the function $f_r(t,n)$. We have determined exact values only in some very special cases; mostly we have upper and lower bounds which are far apart.

In \S2 we establish some general facts about the function $f_r(t,n)$, most notably an asymptotic lower bound: for fixed integers $r,t\geq2$ and any $\varepsilon>0$ we show that $f_r(t,n)>n^{r-1-\varepsilon}$ for all sufficiently large $n$ (Theorem \ref{T7}).

In \S3 we find upper bounds for $f_3(2,n)$, the minimum number of maximal cliques in a triple system of order $n$ and its complement. The main result here is that $f_3(2,n)\leq\lf\frac{n^2}4\rf+5$. For $n\leq12$ this is definitely not the best possible, e.g., $f_3(2,7)=14$ is attained by the Fano plane (Theorem \ref{T11}).

In \S4 we find upper bounds and a few exact values for $f_2(t,n)$ when $t>2$. We show that, if there is a projective plane of order $q$, then $f_2(q+1,n)=q^2+q-1$ for $(q-1)^2<n\leq(q-1)q$, and $f_2(q+1,n)=q^2+q$ for $(q-1)q<n\leq q^2$ (Theorem \ref{T19}). For $q\geq2$ we observe that $f_2(q+1,q^2)=q^2+q$ if and only if a projective plane of order $q$ exists (Theorem \ref{T20}). We obtain an upper bound of the form $f_2(t,n)\leq n+C_t$ whenever a projective plane of order $t$ or $t-1$ exists (Theorem \ref{T21}); in particular see Theorem \ref{T23} for the case $t=3$, Theorem \ref{T25} for $t=4$, Theorems \ref{T26} and \ref{T27} for $t=5$. We have no counterexample to the conjecture that $f_2(3,n)=n+2$ for $n\equiv1\pmod3$ and $f_2(3,n)=n+3$ otherwise (Question \ref{Q24}). 

In \S5 we characterize the graphs $G$ with $c(G)+\overline c(G)=n(G)+2$ (Theorem \ref{T34}(b)), thereby showing that $c(G)+\overline c(G)\geq n(G)+3$ for all graphs outside of a well-defined class of exceptions. In particular, we show that $c(G)+\overline c(G)\geq n(G)+3$ if $\omega(G)+\alpha(G)<n(G)$ (Lemma \ref{L33}, Corollary \ref{C35}(c)).

Some of our results (Theorems \ref{T6}, \ref{T7}, \ref{T8}, \ref{T17}, \ref{T20}, and \ref{T23}, perhaps in a less general formulation) were stated without proof in the earlier paper \cite{GK71} or in the abstract \cite{EGK77}. Some were presented by the third author in invited addresses to the Mid-Atlantic Mathematical Logic Seminar, Dartmouth College, Hanover, New Hampshire, October 2003, and to the CombinaTexas Conference, Houston, Texas, April 2009.

\section{Generalities}

In this section we establish some basic properties of the function $f_r(t,n)$ such as monotonicity (Theorem \ref{T3}), the trivial bounds (Theorem \ref{T6}), and an asymptotic lower bound (Theorem \ref{T7}). We begin by disposing of the trivial cases where $n\leq r+1$. Recall that the maximum size of a $t$-partite graph of order $n$ is attained by the Turán graph $T_{n,t}$.

\begin{theorem}
\label{T2}
For any integers $r,t\geq2$ we have
$$f_r(t,n)=\begin{cases}
t & \text{ if }1\leq n<r,\\
(t-1)n+1 & \text{ if }n=r,\\
\binom{n+1}{2}+(t-2)\binom{n}{2}-e(T_{n,t}) & \text{ if }n=r+1,\\
\lf\frac{(n+1)^2}{4}\rf & \text{ if }n=r+1\text{ and }t=2.\\
\end{cases}$$
\end{theorem}

\begin{proof}

Only the case $n=r+1$ is in need of a proof. Consider a factorization $H_1,\dots,H_t$ of $K_n^{n-1}=\left(V,\binom V{n-1}\right)$. Let $V_i=\{v\in V:V\setminus\{v\}\in E(H_i)\}$ and let $n_i=|V_i|=e(H_i)$, so that $\sum_{i=1}^tn_i=|V|=n$.

First suppose $n_i < n$ for all $i$. Then for $x\in V$ the set $V\setminus\{x\}$ is a maximal clique only for the unique factor $H_i$ having $V\setminus\{x\}$ as an edge, while for $\{x,y\}\in\binom V2$ the set $V\setminus\{x,y\}$ is a maximal clique of $H_i$ just in case $x\notin V_i$ and $y\notin V_i$. Hence the factors $H_i$ have a total of $n$ maximal cliques of size $n-1$ and $(t-2)\binom n 2+\sum_{i=1}^t\binom{n_i}2$ maximal cliques of size $n-2$, so that in this case
\begin{equation}
\label{T2 eqn}
\begin{split}
  \sum_{i=1}^tc(H_i) &= n+(t-2)\binom n2+\sum_{i=1}^t\binom{n_i}2 \\
  &=\binom{n+1}2+(t-2)\binom n2-\left[\binom n2-\sum_{i=1}^t\binom{n_i}2\right]. 
\end{split}
\end{equation}

On the other hand, if $n_1=n$, then
$$\sum_{i=1}^t c(H_i) = 1 + (t-1)\binom n2 = n + (t-2)\binom n2 + \binom{n-1}2,$$
the same number of maximal cliques as when $n_1=n-1$. Hence we may assume that $n_i<n$ for all $i$. The quantity (\ref{T2 eqn}) is minimized when $\binom n2-\sum_{i=1}^t\binom{n_i}2$ is maximized, i.e., when $\binom n2-\sum_{i=1}^t\binom{n_i}2=e(T_{n,t})$. This shows that $f_{n-1}(t,n)=\binom{n+1}2+(t-2)\binom n2-e(T_{n,t})$. Finally, setting $t=2$, we have
$$f_{n-1}(2,n)=\binom{n+1}2-e(T_{n,2})=\binom{n+1}2-\lf\frac{n^2}4\rf=\lf\frac{(n+1)^2}4\rf.$$
\end{proof}

\begin{theorem}
\label{T3}
$f_r(t, n) \leq f_r(t, n+1)$.
\end{theorem}

\begin{proof}
If $H$ is an induced subhypergraph of $G$ then $c(H) \leq c(G)$.
\end{proof}

The inequality is not necessarily strict, e.g., if $n$ is odd then $f_2(n, n) = f_2(n, n+1) = \binom{n+1}{2}$ by Corollary \ref{C15}(b). See Theorem \ref{T19} for more examples.

\begin{question}
\label{Q4}
Does the inequality $f_r(t, n) \leq f_r(t+1, n)$ hold for all integers $r, t \geq 2$ and $n \geq 1$?
\end{question}

At any rate we don't have strict inequality, since $f_2(3, 9) = f_2(4, 9) = 12$ by Theorems \ref{T23} and \ref{T25}.

An obvious way to get a $t$-coloring from a $(t+1)$-coloring is by merging two colors. This can increase the number of maximal monochromatic cliques, i.e., if $G_1$ and $G_2$ are edge-disjoint graphs on the same vertex set $V$, we can have $c(G_1 \cup G_2) > c(G_1) + c(G_2)$. For example, let $V$ be the set of vertices of a unit cube, and let 
\begin{align*}
E(G_1) &= \{\{x,y\} \in \binom{V}{2}:d(x,y) = 1\},\\
E(G_2) &= \{\{x,y\} \in \binom{V}{2}:d(x,y) = \sqrt{2}\}.
\end{align*}
Then $G_1=Q_3$, $G_2=2K_4$, and $G_1 \cup G_2 = K_{2,2,2,2}$, so that $c(G_1)=12$, $c(G_2)=2$, and $c(G_1 \cup G_2) = 16$.

\begin{lemma}
\label{L5}
Consider a coloring of $K_n^r$. For each color $i$, let $c_i$ denote the number of maximal $i$-cliques. Then, for any two colors $i \neq j$, we have $c_ic_j \geq \binom{n}{r-1}$.
\end{lemma}

\begin{proof}
Let $V$ be the vertex set. For each color $i$ let $C_i$ be the set of all maximal $i$-cliques. Each set $X \in \binom{V}{r-1}$ is an $i$-clique; let $f_i(X)$ be some maximal $i$-clique containing $X$. If $i \neq j$, then the map $X \mapsto (f_i(X), f_j(X))$ is an injection from $\binom{V}{r-1}$ to $C_i \times C_j$.
\end{proof}

\begin{theorem}
\label{T6}
For all integers $r, t \geq 2$ and $n \geq r-1$, 

$t\binom{n}{r-1}^\frac{1}{2}\leq f_r(t,n) \leq (t-1) \binom{n}{r-1} + 1$.
\end{theorem}

\begin{proof}
For the upper bound, give all edges the same color.

For the lower bound, let $K_n^r$ be $t$-colored so that $f_r(t, n) = c_1 + \dots + c_t$ where $c_i$ is the number of maximal $i$-cliques. By Lemma \ref{L5} we have $c_ic_j \geq \binom{n}{r-1}$ for $i \neq j$. Hence:

$$(c_1 c_2)(c_2 c_3) \dots (c_t c_1) \geq \binom{n}{r-1}^t;$$
$$c_1 c_2 \dots c_t \geq \binom{n}{r-1}^{t/2};$$
$$f_r(t,n) = c_1 + c_2 + \dots + c_t \geq t(c_1 c_2 \dots c_t)^{1/t} \geq t \binom{n}{r-1}^{1/2}.$$
\end{proof}

\begin{theorem}
\label{T7}
For any integers $s \geq 1$ and $t \geq 2$, and for any $\varepsilon > 0$, we have $f_{s+1}(t, n) > n^{s - \varepsilon}$ for all sufficiently large $n$.
\end{theorem}

\begin{proof}
Choose $\delta$ so that $0<\delta<\min\{\varepsilon,s\}$. Let $\alpha=\frac{s-\delta}s$ and $k=\left\lceil\frac{s-\delta}{1-\alpha}\right\rceil$. Define $M$ so that $$m\geq M\implies\frac{m(m-1)(m-2)\cdots(m-s+1)}{m^s}\geq\frac{t-\frac32}{t-1}.$$
  By Ramsey's theorem there is an integer $R\geq k$ such that any $t$-coloring of the edges of the complete hypergraph $K_R^{s+1}$ contains a monochromatic clique of size $k$. Let $\beta=\min\left\{\frac{t-\frac32}{s!},\frac1{\binom Rk}\right\}$.

Let $n\in\mathbb N$ and consider a $t$-coloring of the edges of $K_n^{s+1}$. Let $c_i$ be the number of maximal $i$-cliques. Let $m$ be the size of the largest monochromatic clique.

\item[\emph{\textbf{Case 1.} $m\geq n^\alpha$.}] 

If there is a $j$-clique of size $m$, then for each color $i\neq j$ there are at least $\binom ms$ maximal $i$-cliques. If $n^\alpha\geq M$ then $$\sum_{i=1}^tc_i>(t-1)\binom ms\geq\frac{t-\frac32}{s!}m^s\geq\beta m^s\geq\beta n^{\alpha s}=\beta n^{s-\delta}.$$

\item[\emph{\textbf{Case 2.} $m\leq n^\alpha$.}] 

Let $p$ be the number of monochromatic cliques of size $k$. Since every monochromatic clique is contained in a maximal monochromatic clique, $p\leq\binom mk\sum_{i=1}^tc_i$. On the other hand, if $n\geq R$, then $p\geq\binom nk/\binom Rk\geq\beta\binom nk$ and we have $$\sum_{i=1}^tc_i\geq\frac{\beta\binom nk}{\binom mk}\geq\frac{\beta n^k}{m^k}\geq\frac{\beta n^k}{n^{\alpha k}}\geq\beta n^{s-\delta}.$$

Thus $f_{s+1}(t,n)\geq\beta n^{s-\delta}>n^{s-\varepsilon}$ for sufficiently large $n$.
\end{proof}

\section{A triple system and its complement}

In this section we find upper bounds for $f_3(2,n)$. For $n>12$ our best
result is $f_3(2,n)\leq\lf\frac{n^2}4\rf+5$ (Theorem
\ref{T12}). We know the exact value of $f_3(2,n)$ only for $n\leq7$. By tedious
case analysis (which we omit) we have verified that Theorem \ref{T8} is optimal
for $2\leq n\leq6$, so that
$f_3(2,n)=\lf\frac{(n+1)^2}4\rf$ in those cases, and
that Theorem \ref{T11} is optimal for $n=7$, so that $f_3(2,7)=14$. See Theorem
\ref{T2} for the trivial cases $n\leq4$.

\begin{theorem}
\label{T8}
$f_3(2, n) \leq \lf\frac{(n+1)^2}{4}\rf$ for $n \geq 2$.
\end{theorem}

\begin{proof}
Let $n = n_1 + n_2$ where $n_1, n_2 > 0$. Let $V = V_1 \cup V_2, V_1 \cap V_2 = \varnothing, |V_i| = n_i$, and consider the triple system $H = (V, E)$ where $E = \{e \in \binom{V}{3} : |e \cap V_1| \geq 2\}$. Then $c(H) = n_2 + \binom{n_2}{2}$; the maximal cliques are the sets $V_1 \cup \{v\}$, $v \in V_2$, and the $2$-element subsets of $V_2$. Taking $n_1 = \lf\frac{n}{2}\rf$ and $n_2 = \lc\frac{n}{2}\rc$, we see that
$$f_3(2, n) \leq c(H) + \overline{c}(H) = \lc\frac{n}{2}\rc  + \binom{\lc\frac{n}{2}\rc}{2} + \lf\frac{n}{2}\rf  + \binom{\lf\frac{n}{2}\rf}{2} = \lf\frac{(n+1)^2}{4}\rf.$$

Alternatively, let $H=(V,E)$ where $V=[n]$ and $E=\{\{x,y,z\}\in\binom V3:x<y<z,\ y\text{ odd}\}$. For integers $a < b$ of the same parity, let $X_{a,b}$ be the set consisting of $a$ and $b$ and all integers of the
opposite parity between $a$ and $b$; e.g., $X_{3,9}=\{3,4,6,8,9\}$. Then
the maximal cliques (anticliques) of $H$ are just the sets $X_{a,b}\cap[n]$ where $a$ and $b$ are even (odd) integers, $0\leq a < b \leq n+1$. It follows that
$$c(H)+\overline{c}(H) = \binom{\left\lceil\frac
n2\right\rceil+1}2+\binom{\left\lfloor\frac
n2\right\rfloor+1}2=\left\lfloor\frac{(n+1)^2}4\right\rfloor.$$
\end{proof}

In fact $f_3(2, n) = \lf\frac{(n+1)^2}{4}\rf$ holds for $2 \leq n \leq 6$ but not for $n \geq 7$; see Theorems \ref{T11} and \ref{T12}.

\begin{lemma}
\label{L9}
A triple system $H$ of order $n \geq 2$ can be extended to a triple system $H^*$ of order $n+1$ such that $c(H^*) = c(H) + d(H) + 1$, $d(H^*) = d(H) + 1$, $\overline{c}(H^*) = \overline{c}(H)$, and $\overline{d}(H^*) = \overline{d}(H)$.
\end{lemma}

\begin{proof}
Let $H = (V, E)$ be a triple system of order $n \geq 2$. Choose a vertex $v \in V$ with $D(v) = d(H)$. Choose a new vertex $v^* \notin V$ and let $H^* = (V^*, E^*)$ where $V^* = V \cup \{v^*\}$ and $E^* = E \cup \{\{x, y, v^*\} : \{x, y, v\} \in E\}$.

The maximal anticliques of $H^*$ are $(1)$ the maximal anticliques of $H$ not containing $v$ and $(2)$ the sets $X \cup \{v^*\}$ where $X$ is a maximal anticlique of H containing $v$. Plainly $\overline{c}(H^*) = \overline{c}(H)$ and $\overline{d}(H^*) = \overline{d}(H)$.

The maximal cliques of $H^*$ are $(1)$ the maximal cliques of $H$, and $(2)$ the sets $(X \setminus \{v\}) \cup \{v^*\}$ where $X$ is a maximal clique of $H$ containing $v$, and $(3)$ the set $\{v, v^*\}$. Now it can be seen that $c(H^*) = c(H) + D_H(v) + 1 = c(H) + d(H) + 1$, and that $D_{H^*}(v) = D_{H^*}(v^*) = D_H(v)+1 = d(H)+1$, while $D_{H^*}(u) \geq D_H(u) + 1 \geq d(H) + 1$ for $u \in V \setminus \{v\}$, whence $d(H^*) = d(H) + 1$.
\end{proof}

\begin{lemma}
\label{L10}
If there is a triple system $H$ of order $n \geq 2$ with $c(H) + \overline{c}(H) = k$ and $d(H) = \lf\frac{m}{2}\rf$ and $\overline{d}(H) = \lc\frac{m}{2}\rc$, then there is a triple system $H'$ of order $n+1$ with $c(H') + \overline{c}(H') = k + \lf\frac{m+2}{2}\rf$ and $d(H') = \lf\frac{m+1}{2}\rf$ and $\overline{d}(H') = \lc\frac{m+1}{2}\rc$.
\end{lemma}

\begin{proof}
By Lemma \ref{L9} there is a triple system $H^*$ of order $n+1$ such that $c(H^*) + \overline{c}(H^*) = c(H) + d(H) + 1 + \overline{c}(H) = k + \lf\frac{m+2}{2}\rf$ and $d(H^*) = d(H)+1 = \lf\frac{m+2}{2}\rf = \lc\frac{m+1}{2}\rc$ and $\overline{d}(H^*) = \overline{d}(H) = \lc\frac{m}{2}\rc = \lf\frac{m+1}{2}\rf$. Let $H'$ be the complement of $H^*$.
\end{proof}

\begin{theorem}
\label{T11}
$f_3(2, n) \leq \lf\frac{(n+1)^2}{4}\rf - 2$ for $n \geq 7$.
\end{theorem}

\begin{proof}
We prove by induction that for each $n \geq 7$ there is a triple system $H$ of order $n$ with $c(H) + \overline{c}(H) = \lf\frac{(n+1)^2}{4}\rf - 2$ and $d(H) = \lf \frac{n}{2} \rf$ and $\overline{d}(H) = \lc \frac{n}{2} \rc$. The inductive step follows from Lemma \ref{L10} since

$$\lf\frac{(n+1)^2}{4}\rf - 2 + \lf \frac{(n+2)}{2} \rf = \lf \frac{(n+2)^2}{4} \rf - 2.$$

For the base case $n=7$, consider the Fano plane as a hypergraph H with the points as vertices and the lines as edges. Then the maximal cliques are the lines, and the maximal anticliques are the complements of lines. Thus $c(H) = \overline{c}(H) = 7$, and $d(H) = 3$ and $\overline{d}(H) = 4$; each vertex is in exactly three maximal cliques and four maximal anticliques.
\end{proof}

In fact $f_3(2,7) = 14$; the proof of $f_3(2,7) \geq 14$ is a tedious case analysis which we omit. We do not know if equality holds in Theorem \ref{T11} when $8 \leq n \leq 14$. For $n \geq 15$ a better bound is given by Theorem \ref{T12}.

\begin{theorem}
\label{T12}
$f_3(2, n) \leq \lf\frac{n^2}{4}\rf + 5$.
\end{theorem}

\begin{proof}
The cases with $n \leq 5$, while easily verified, are of no interest. For $n \geq 6$ we prove by induction that there is a triple system $H$ of order $n$ with $c(H) + \overline{c}(H) = \lf \frac{n^2}{4} \rf + 5$ and $d(H) = \lf \frac{n-1}{2} \rf$ and $\overline{d}(H) = \lc \frac{n-1}{2} \rc$. The inductive step follows from Lemma \ref{L10} since 
$$\lf \frac{n^2}{4} \rf + 5 + \lf \frac{n+1}{2} \rf = \lf \frac{(n+1)^2}{4} \rf + 5.$$

Now consider $H = (V, E)$ where $V=[6]$ and $E = \{\{1,3,6\},\\ \{1, 4, 5\}, \{2,3,5\}, \{2,3,6\}, \{2,4,5\}, \{2,4,6\}\}$. 

We can visualize this as the $6$ vertices and $6$ of the $8$ faces of a regular octahedron, the two missing faces meeting at one point. There are $9$ maximal cliques, namely, the $6$ edges and the $3$ antipodal pairs 
$$\{1, 2\}, \{3, 4\}, \{5, 6\};$$ and there are $5$ maximal anticliques, namely, the missing faces $$\{1,3,5\}, \{1,4,6\},$$ and the complements of the antipodal pairs, $$\{3,4,5,6\}, \{1,2,5,6\}, \{1,2,3,4\}.$$ Thus 
$$c(H) + \overline{c}(H) = 9+5 = 14,$$ 
while $$d(H) = D(1) = 3$$ and $$\overline{d}(H) = \overline{D}(2) = 2.$$ The complement of $H$ satisfies the requirements of the base case $n=6$. 
\end{proof}

\section{Edge-colored graphs and projective planes}

The main results of this section are upper bounds for $f_2(t,n)$ from recursive constructions using projective planes. For example, we use projective planes of orders $4$ and $5$ to show that $f_2(5,n)\leq n+7$ for $n\geq37$ (Theorem \ref{T27}). From Theorem \ref{T20} we see that determining the exact value of $f_2(t,n)$ in all cases will be at least as hard as determining all possible orders of finite projective planes.

\begin{theorem}
\label{T13}
For all integers $t \geq 2$ and $n \geq 1$,

$t \sqrt{n} \leq (t-2)\lc\sqrt{n}\rc + \lc\sqrt{4n} \hspace{0.1cm} \rc \leq f_2(t, n) \leq (t-1)n + 1$.
\end{theorem}

\begin{proof}
The upper bound is from Theorem \ref{T6}. For the lower bound, let $f_2(t, n) = c_1 + c_2 + \dots + c_t$ where $c_i$ is the number of maximal $i$-cliques in some $t$-coloring of $K_n$; we may assume that $c_1 \leq c_2 \leq \dots \leq c_t$. By Lemma \ref{L5} we have $c_1c_2 \geq n$. It follows that $c_i \geq c_2 \geq \lc \sqrt{n} \hspace{0.1cm} \rc$ for $i = 3, \dots, t$, and $c_1 + c_2 \geq  \lc 2 \sqrt{c_1 c_2} \hspace{0.1cm} \rc \geq \lc \sqrt{4n} \hspace{0.1cm} \rc$, whence $f_2(t, n) = c_1 + c_2 + \dots + c_t \geq (t-2) \lc \sqrt{n} \hspace{0.1cm} \rc + \lc \sqrt{4n} \hspace{0.1cm} \rc \geq t \sqrt{n}$.
\end{proof}

\begin{theorem}
\label{T14}
$f_2(t, n) \geq tn - \binom{n}{2}$, with equality if and only if $n \leq 2\lc\frac{t}{2}\rc$.
\end{theorem}

\begin{proof}
Let the edges of $K_n$ be colored with $t$ colors. Consider the graph $G = G_1 \cup G_2 \cup \dots \cup G_t$ where $G_1, \dots, G_t$ are vertex-disjoint graphs such that $G_i$ is isomorphic to the spanning subgraph of $K_n$ formed by the edges of color $i$. Let $k$ be the number of components of $G$. Then
\begin{equation}
\label{T14 eqn}
 \sum_{i=1}^{t} c(G_i) = c(G) \geq k \geq |V(G)| - |E(G)| = tn - \binom{n}{2}. 
\end{equation}
For equality to hold in (\ref{T14 eqn}) we must have $c(G) = k$, meaning that each component of $G$ is a clique, and also $k = |V(G)| - |E(G)|$, meaning that $G$ is acyclic. So equality holds just in case each component of $G$ is $K_1$ or $K_2$, which means that the given coloring is a proper edge-coloring of $K_n$. Of course, $K_n$ is $t$-edge-colorable if and only if $n \leq 2 \lc \frac{t}{2} \rc$.
\end{proof}

\begin{corollary}  
\label{C15}    
  \[
\begin{array}{ll}
 (a) & f_2(n+1, n) = \frac{n(n+3)}{2}$.$\\
 (b) & $ If $n$ is odd, then $ f_2(n, n+1) = \binom{n+1}{2} = f_2(n, n)$.$\\
 (c) & $ If $n$ is even, then $ f_2(n, n+1) > \binom{n+1}{2} = f_2(n, n)$.$
\end{array}
  \]
\end{corollary}

\begin{lemma}
\label{L16}
Let $G$ be a $t$-colored complete graph of order $m$, and let $c = c_1 + \cdots + c_t$ where $c_i$ is the number of maximal $i$-cliques in $G$. If some vertex $v$ of $G$ is in a unique maximal $i$-clique for each color $i$, then $f_2(t, n+m-1) \leq f_2(t, n) + c - t$ for all $n \in \mathbb{N}$.
\end{lemma}

\begin{proof}
Let $H$ be an optimally $t$-colored complete graph of order $n$, so that $f_2(t, n) = c_1^H + \dots + c_t^H$ where $c_i^H$ is the number of maximal $i$-cliques in $H$. Replace the vertex $v$ of $G$ with a copy of $H$. In the resulting $t$-colored complete graph of order $n+m-1$, the number of maximal $i$-cliques meeting $V(H)$ is $c_i^H$, while the number of maximal $i$-cliques disjoint from $V(H)$ is $c_i - 1$.
\end{proof}

\begin{theorem}
\label{T17}
If there is a projective plane of order $q$, then $f_2(q+1, n+q^2-1) \leq f_2(q+1, n) + q^2 - 1$ for all $n \in \mathbb{N}$.
\end{theorem}

\begin{proof}
Let $P$ be the point set of a projective plane of order $q$. Choose a line $l_0$ and let $x_0, x_1, \dots, x_q$ be the points on $l_0$. Let $G$ be the complete graph of order $q^2$ with vertex set $P \setminus \{x_0, x_1, \dots, x_q\}$. Color the edges of $G$ with colors $0, 1, \dots, q$ by assigning color $i$ to an edge $uv$ if the points $u, v, x_i$ are collinear. For each color $i$ there are $q$ maximal $i$-cliques corresponding to the lines other than $l_0$ through the point $x_i$, and each vertex of $G$ is in just one of them. We obtain the desired inequality by setting $t = q + 1$, $m = q^2$, $c = q^2 + q$ in Lemma \ref{L16}.
\end{proof}

The following unpublished observation by Graham and Van Lint is included by permission.

\begin{theorem}
\label{T18}
(Ronald L. Graham and Jack van Lint)
If there is a projective plane of order $q$, then $f_2(q, n+q^2) \leq f_2(q, n) + q^2$ for all $n \in \mathbb{N}$.
\end{theorem}

\begin{proof}
Let $P$ be the point set of a projective plane of order $q$. Choose a line $l_0$; let $x_0, x_1, \dots, x_q$ be the points on $l_0$, and let $l_1, \dots, l_q$ be the other lines through $x_0$. Let $G$ be the complete graph of order $q^2 + 1$ with vertex set $P \setminus \{x_1, \dots, x_q\}$. For any edge $uv$ of $G$, let $l(uv)$ be the line through $u$ and $v$. Color the edges of $G$ with colors $1, \dots, q$ by assigning color $i$ to an edge $uv$ if either $l(uv) = l_i$ or else $l(uv)$ meets $l_0$ at the point $x_i$. For each color $i$ there are $q+1$ maximal $i$-cliques; $q$ of them correspond to the lines other than $l_0$ through $x_i$; the remaining one, corresponding to the line $l_i$, is the only maximal $i$-clique containing $x_0$. We obtain the desired inequality by setting $t=q$, $m=q^2+1$, $c=q^2+q$, and $v=x_0$ in Lemma \ref{L16}.
\end{proof}

\begin{theorem}
\label{T19}
If there is a projective plane of order $q$, then:
  \[
\begin{array}{ll}
 (a) & f_2(q+1, n) = q^2 + q - 1$ whenever $(q-1)^2 < n \leq (q-1)q$;$\\
 (b) & f_2(q+1, n) = q^2 + q$ whenever $(q-1)q < n \leq q^2$.$
\end{array}
  \]
\end{theorem}

\begin{proof}
If $n > (q-1)^2$ then by Theorem \ref{T13} we have 
\begin{gather}
f_2(q+1, n) \geq (q-1) \lc \sqrt{(q-1)^2 + 1} \rc + \lc \sqrt{4(q-1)^2 + 4} \rc \notag \\
= (q-1)q + (2q-1) = q^2 + q - 1; \notag
\end{gather}
moreover, if $n > (q-1)q$, then
\begin{gather}
f_2(q+1, n) \geq (q-1) \lc \sqrt{(q-1)q + 1} \rc + \lc \sqrt{4(q-1)q + 4} \rc \\
= (q-1)q + 2q = q^2 + q.
\end{gather}

On the other hand, if there is a projective plane of order $q$, and if $n \leq q^2$, then
$$f_2(q+1, n) \leq f_2(q+1, q^2) \leq f_2(q+1, 1) + q^2 - 1 = q^2 + q$$
by Theorem \ref{T17}. This completes the proof of $(b)$.

Take a projective plane of order $q$. Choose a line $l_0$ with points\\
$x_0, x_1, \dots, x_q$ and another line $l_1$ through $x_0$. Let $G$ be the complete graph of order $q^2-q$ whose vertices are the points not on $l_0$ or $l_1$. Assign the color $i$ to an edge $uv$ of $G$ if the line through $u$ and $v$ meets $l_0$ at $x_i$. Now maximal $i$-cliques correspond to lines through $x_i$ other than $l_0$ and $l_1$, so the number of maximal $i$-cliques is $q-1$ if $i=0$ and $q$ if $i \in \{1, \dots, q\}$, for a total of $q^2+q-1$. Hence, for $n \leq q^2-q$, we have 
$$f_2(q+1, n) \leq f_2(q+1, q^2-q) \leq q^2 + q - 1.$$ 
This completes the proof of $(a)$.
\end{proof}

\begin{theorem}
\label{T20}
For $n \geq 2$ the following statements are equivalent:
  \[
\begin{array}{ll}
 (a) & f_2(n+1, n^2) = n^2 + n$;$\\
 (b) & f_2(n+1, n^2) \leq n^2 + n$;$\\
 (c) & $there is a projective plane of order $n$.$
\end{array}
  \]
\end{theorem}

\begin{proof}
We have $(c) \implies (a)$ by Theorem \ref{T19}$(b)$, and $(a) \implies (b)$ is trivial; we have to show $(b) \implies (c)$.

Assume that $n \geq 2$ and $f_2(n+1, n^2) \leq n^2 + n$. Consider a factorization $K_{n^2} = G_0 \cup G_1 \cup \dots \cup G_n$ with $c_0 + c_1 + \dots + c_n \leq n^2 + n$, where $c_i = c(G_i)$. Since $c_ic_j \geq n^2$ for $i \neq j$, the geometric mean of $c_0, c_1, \dots, c_n$ is at least $n$, that is, 
$$n \leq (c_0 c_1 \dots c_n)^{\frac{1}{n+1}} \leq \frac{c_0 + c_1 + \dots + c_n}{n+1} \leq n,$$
whence $c_0 = c_1 = \dots = c_n = n$.

For each $i$, there are exactly $n$ maximal $i$-cliques, which are pairwise disjoint and contain $n$ vertices each. If $X$ is a maximal $i$-clique and $Y$ a maximal $j$-clique, $i \neq j$, then $|X \cap Y| = 1$. Hence there is a projective plane of order $n$; the points are the vertices of $K_{n^2}$ and the numbers $0, 1, \dots, n$; the lines are the set $\{0, 1, \dots, n\}$ and the sets $X \cup \{i\}$ where $X$ is a maximal $i$-clique.
\end{proof}

In view of Theorem \ref{T20}, a computation by Lam, Thiel, and Swiercz \cite{LTS89} shows that $f_2(11,100)>110$. On the other hand, $$f_2(11,100)\leq f_2(11,122)\leq f_2(11,1)+121=132$$ by Theorem \ref{T18}.

\begin{theorem}
\label{T21}
(Due to Ronald L. Graham and Jack van Lint in the case of a projective plane of order $t$.) If there is a projective plane of order $t$ or $t-1$, then there is a constant $C_t$ such that $f_2(t,n)\leq n+C_t$ for all $n\in\mathbb N$.
\end{theorem}

\begin{proof}
If there is a projective plane of order $t-1$, the result follows from Theorem \ref{T17} with
$$C_t = \max\{f_2(t, n) - n: 1 \leq n \leq t^2 - 2t\};$$
if there is a projective plane of order $t$, it follows from Theorem \ref{T18} with
$$C_t = \max\{f_2(t, n) - n: 1 \leq n \leq t^2\}.$$
\end{proof}

\begin{question}
\label{Q22}
For each integer $t \geq 2$, is there a constant $C_t$ such that $f_2(t, n) \leq n+C_t$ for all $n \in \mathbb{N}$?
\end{question}

If $C_t$ exists then $C_t \geq \binom{t}{2}$ since, when $n = t-1$, we have $f_2(t, n) = f_2(n+1, n) = \frac{n(n+3)}{2} = n + \binom{n+1}{2} = n + \binom{t}{2}$ by Corollary \ref{C15}$(a)$. In fact the inequality $f_2(t, n) \leq n + \binom{t}{2}$ holds for $2 \leq t \leq 5$ (Theorems \ref{T1}, \ref{T23}, \ref{T25}, \ref{T26}). On the other hand it fails for $t=6$ and $t=8$, since $f_2(6,7) = 23$ and $f_2(8, 9) = 38$; we omit the details.

Of course a negative answer to Question \ref{Q22} would imply a negative answer to Question \ref{Q4} as well.

\begin{theorem}
\label{T23}
  \[
    f_2(3, n) \leq\left\{
                \begin{array}{ll}
                  n+2 & $ if $n\equiv1\pmod3$,$\\
                  n+3 & $ otherwise,$
                \end{array}
              \right.
  \]
with equality at least for $n \leq 10$.
\end{theorem}

\begin{proof}
By Theorem \ref{T17} with $q=2$, to verify the inequality it suffices to observe that it holds for $n \leq 3$. To verify equality for all $n \leq 10$, it suffices to show that $f_2(3, 8) \geq 11$ and $f_2(3, 9) \geq 12$; we omit the details.
\end{proof}

\begin{question}
\label{Q24}
Does equality hold in Theorem \ref{T23} for all $n$?
\end{question}

\begin{theorem}
\label{T25}
  \[
    f_2(4, n) \leq\left\{
                \begin{array}{ll}
                  n+3 & $if $n\equiv 1 \pmod 8$,$\\
                  n+4 & $if $n\equiv 0 \pmod 8$,$\\
                  n+5 & $if $n\equiv 2, 6, 7 \pmod 8$,$\\
                  n+6 & $if $n\equiv 3, 4, 5 \pmod 8$,$\\
                \end{array}
              \right.
  \]
  with equality at least for $n \leq 10$.
\end{theorem}

\begin{proof}
Use $q=3$ in Theorem \ref{T17}.
\end{proof}

\begin{theorem}
\label{T26}
$f_2(5, n) \leq n+10$ for all $n$.
\end{theorem}

\begin{proof}
By Theorem \ref{T17} with $q=4$ we have $f_2(5, n+15) \leq f_2(5, n) + 15$; hence it will suffice to prove the inequality $f_2(5, n) \leq n + 10$ for $ n \leq 15$. By Theorem \ref{T14} we have $f_2(5, n) = 5n - \binom{n}{2}$ for $n \leq 6$, and by Theorem \ref{T19} with $q=4$ we have $f_2(5, n) = 19$ for $10 \leq n \leq 12$ and $f_2(5, n) = 20$ for $13 \leq n \leq 16$. We leave it as an exercise for the reader to verify that $f_2(5, 7) \leq 17$ and $f_2(5, 8) \leq 18$. (In fact it is easy to show that $f_2(5, 7) = 17$ and $f_2(5, 8) = f_2(5, 9) = 18$.)
\end{proof}

\begin{theorem}
\label{T27}
$f_2(5,n) \leq n+7$ for $n \geq 37$. In fact, if $n \geq 37$, then
  \[
\begin{array}{ll}
 (a) & f_2(5, n) \leq n + 4 $ if $ n \equiv 1 \pmod 5$,$\\
 (b) & f_2(5, n) \leq n + 7 $ if $ n \equiv 2 \pmod 5$,$\\
 (c) & f_2(5, n) \leq n + 7 $ if $ n \equiv 3 \pmod 5$,$\\
 (d) & f_2(5, n) \leq n + 6 $ if $ n \equiv 4 \pmod 5$,$\\
 (e) & f_2(5, n) \leq n + 5 $ if $ n \equiv 0 \pmod 5$.$\\
\end{array}
  \]
\end{theorem}

\begin{proof}
Only $(a)$ and $(b)$ need proof, as $(c)-(e)$ follow directly from $(a)$. We will use the facts that $f_2(5,1)=5$ and $f_2(5,2)=9$ by Theorem \ref{T2} (or Theorem \ref{T14}), and $f_2(5,12)=19$ by Theorem \ref{T19}$(a)$ with $q=4$.

We have $f_2(5, n+15) \leq f_2(5, n) + 15$ by Theorem \ref{T17} with $q=4$, and $f_2(5, n+25) \leq f_2(5, n)+25$ by Theorem \ref{T18} with $q=5$. It follows that $f_2(5, n+5k) \leq f_2(5, n) + 5k$ for every integer $k \geq 8$.

If $n \equiv 1 \pmod 5$ and $n > 37$, then $n = 1 + 5k$ for some integer $k \geq 8$, and
$$f_2(5, n) = f_2(5, 1+5k) \leq f_2(5,1)+5k = 5+5k = n+4.$$

If $n \equiv 2 \pmod 5$ and $n > 37$, then $n = 2 + 5k$ for some integer $k \geq 8$, and
$$f_2(5,n) = f_2(5, 2+5k) \leq f_2(5,2)+5k = 9+5k = n+7.$$

Finally, if $n = 37$ then $n \equiv 2 \pmod 5$ and
$$f_2(5,n) = f_2(5, 12+25) \leq f_2(5, 12) + 25 = 19+25 = n+7.$$
\end{proof}

\section{When $c(G)+\overline c(G)=n(G)+2$}

For a graph $G$ we define $\tau(G)=c(G)+\overline c(G)-n(G)$. In this section we characterize the graphs $G$ with $\tau(G)=2$. In doing so we also prove the characterization of graphs with $\tau(G)=1$ which was proved more simply in \cite{GK71}. We need the following improved version of the result of \cite{GK71} that $\tau(G)\geq1$ for every graph $G$.

\begin{theorem}
\label{T28}
If $G$ is a graph of order $n$, then 

$c(G) + \overline{c}(G) \geq n + d(G) + \overline{d}(G) - 1 \geq n + 1$.
\end{theorem}

\begin{proof}
We use induction on $n$. Let $G$ be a graph of order $n$; let $c=c(G)$, $\overline c=\overline c(G)$, $d=d(G)$, $\overline d=\overline d(G)$. We may assume that there is a vertex $v$ such that $N(v) \neq \varnothing \neq \overline{N}(v)$; otherwise $G$ or $\overline{G}$ is a complete graph and the result is clear. Let $G_1 = G[N(v)]$ and $G_2 = G[\overline{N}(v)]$; Let $n_i = n(G_i)$, $c_i=c(G_i)$, $\overline{c}_i = \overline{c}(G_i)$, $d_i = d(G_i)$, $\overline{d}_i = \overline{d}(G_i)$. By the inductive hypothesis, $c_i + \overline{c}_i \geq n_i + d_i + \overline{d}_i - 1$.

Choose a vertex $w$ of $G_2$ which is in exactly $d_2$ maximal cliques of $G_2$. The number of maximal cliques of $G$ containing $v$ is equal to $c_1$; the number of maximal cliques containing $w$ is at least $d$; the number of maximal cliques containing neither $v$ nor $w$ is at least $c_2 - d_2$. Since no clique contains both $v$ and $w$,
$c \geq c_1 + d + c_2 - d_2.$ Similarly, $\overline{c} \geq \overline{c}_2 + \overline{d} + \overline{c}_1 - \overline{d}_1$.

Adding these two inequalities we get
\[
\begin{array}{ll}
c + \overline{c} \geq (c_1 + \overline{c}_1) + (c_2 + \overline{c}_2) + d + \overline{d} - d_2 - \overline{d}_1\\
\geq (n_1 + d_1 + \overline{d}_1 - 1) + (n_2 + d_2 + \overline{d}_2 - 1) + d + \overline{d} - d_2 - \overline{d}_1\\
= n_1 + n_2 + (d_1-1) + (\overline{d}_2 - 1) + d + \overline{d}\\
\geq n_1 + n_2 + d + \overline{d} = n + d + \overline{d} - 1.
\end{array}
\]
\end{proof}

We define four classes of graphs, unimaginatively named after their smallest members. As motivation note that $K_1$ is the smallest graph $G$ with $\tau(G)=1$, while $P_4$, $C_4$, and $\overline C_4$ are the smallest graphs with $\tau(G)=2$. Recall that $G$ is a \emph{split graph} if $V(G)$ is the union of a clique and an anticlique. We begin by defining two classes of split graphs.

A graph $G$ is $K_1$\emph{-like} if $V(G)=X\cup Y$ where $X$ is a clique, $Y$ is an anticlique, and $X\cap Y\neq\varnothing$; equivalently, if $\omega(G)+\alpha(G)=n(G)+1$.

A graph $G$ is $P_4$\emph{-like} if $V(G)=X\cup Y$ where $X$ is a maximal clique, $Y$ is a maximal anticlique, and $X\cap Y=\varnothing$.

\begin{lemma}
\label{L29}
If $G$ is a split graph then just one of the following statements holds:
\[
\begin{array}{ll}
(a) & \omega(G)+\alpha(G)=n(G)+1$, $G$ is $K_1$-like, and $\tau(G)=1;\\
(b) & \omega(G)+\alpha(G)=n(G)$, $G$ is $P_4$-like, and $\tau(G)=2.
\end{array}
\]
\end{lemma}

A graph $G$ is $C_4$\emph{-like} if there is a $4$-element set $U=\{p,q,r,s\}\subseteq V(G)$ such that $pq,qr,rs,ps\in E(G)$ while $pr,qs\notin E(G)$, and $V(G)\setminus U=X\cup Y$ where $X$ is a clique and $Y$ is an anticlique, each vertex in $X$ is joined to $p$ and $q$ and at least one more vertex in $U$, and no vertex in $Y$ is joined to any vertex in $U$.

A graph $G$ is $\overline C_4$\emph{-like} if $\overline G$ is $C_4$-like.

\begin{lemma}
\label{L30}
If a graph $G$ is $P_4$-like or $C_4$-like or $\overline{C}_4$-like, then $\tau(G)=2$ and $\omega(G) + \alpha(G) = n(G)$.
\end{lemma}

\begin{lemma}
\label{L31}
Let $G$ be a graph. If $\tau(G) \leq 2$ and $\omega(G) + \alpha(G) \geq n(G)$, then $G$ is $C_4$-like or $\overline{C}_4$-like or a split graph.
\end{lemma}

\begin{proof}
Let $G=(V, E)$ and $n = |V|$. Let $X$ be a clique and $Y$ an anticlique with $|X|+|Y|=n$. We may assume that $X \cup Y \neq V$; let $V \setminus(X \cup Y) = \{s\}$, $X \cap Y = \{p\}$, $X_0 = X \setminus\{p\}$, $Y_0 = Y \setminus \{p\}$. Since the statement of the lemma is invariant under complementation, we may assume that $p$ is joined to $s$. We may assume that some vertex $q \in X_0$ is not joined to $s$, and some vertex $r \in Y_0$ is joined to $s$, as otherwise $G$ is a split graph. Then $D(p) \geq 2$, since $\{p, q\}$ and $\{p, s\}$ are contained in different maximal cliques.
\par
\item[\emph{Case 1. $q$ is not joined to $r$. }]

Then $\overline{D}(q) \geq 2$. From this and $D(p) \geq 2$ we have
\begin{equation}
\label{T31 eqn a}
n+2 \geq c(G) + \overline{c}(G) \geq \sum_{y \in Y}D(y) + \sum_{x \in X}\overline{D}(x) \geq (|Y|+1) + (|X|+1) = n+2.
\end{equation}
Since equality holds in (\ref{T31 eqn a}), we have $D(y)=1$ for each $y \in Y_0$, $\overline{D}(x) = 1$ for each $x \in X \setminus \{q\}$, and every maximal anticlique meets $X$.

If some vertex $x \in X_0 \setminus \{q\}$ is not joined to $s$, then, according as $x$ is joined to $r$ or not, either $D(r) \geq 2$ or $\overline{D}(x) \geq 2$, neither of which is possible. Hence $(X \setminus \{q\}) \cup \{s\}$ is a clique.

Now suppose $q$ is joined to some vertex $y \in Y_0$. If $y$ is joined to $s$, then $D(y) \geq 2$, which is impossible; otherwise $\{y, s\}$ is contained in a maximal anticlique disjoint from $X$, also impossible. Hence $q$ is joined to no vertex in $Y_0$, so $Y_0 \cup \{q\}$ is an anticlique and $G$ is a split graph.

\item[\emph{Case 2. $q$ is joined to $r$.}]

Then $D(r) \geq 2$. From this and $D(p) \geq 2$ we have
\begin{equation}
\label{T31 eqn b}
n+2 \geq c(G) + \overline{c}(G) \geq \sum_{y \in Y}D(y) + \sum_{x \in X}\overline{D}(x) \geq (|Y| + 2) + |X| = n + 2.
\end{equation}
Since equality holds in (\ref{T31 eqn b}), we have $D(y) = 1$ for each $y \in Y_0 \setminus \{r\}$, $\overline{D}(x) = 1$ for each $x \in X$, and every maximal anticlique meets $X$.

Each vertex $x \in X$ is joined to $r$ or to $s$, as otherwise we would have $\overline{D}(x) \geq 2$.

If $s$ were joined to some vertex $y \in Y_0 \setminus \{r\}$, then we would have $D(y) \geq 2$ or $\overline{D}(q) \geq 2$ according as $y$ is joined to $q$ or not. Hence $s$ is not joined to any vertex in $Y_0 \setminus \{r\}$.

Now the anticlique $(Y_0 \setminus \{r\}) \cup \{s\}$ is contained in a maximal anticlique, which must meet $X$ but cannot contain $p$. Hence there is a vertex $q' \in X_0$ which is not joined to $s$ or to any vertex in $Y_0 \setminus \{r\}$. Of course $q'$ is joined to $r$, since $\overline{D}(q') = 1$.

It is easy to check that $G$ is $C_4$-like, with $p, q', r, s, X_0 \setminus \{q'\}, Y_0 \setminus \{r\}$ playing the roles of $p, q, r, s, X, Y$ in the definition.
\end{proof}

\begin{lemma}
\label{L32}
If $G$ is a graph and $H = G[\overline{N}(v)]$, where $v$ is a vertex such that $N(v)$ is a clique and $\overline{N}(v) \neq \varnothing$, then $\tau(H) \leq \tau(G)$. 
\end{lemma}

\begin{proof}
Let $X = N(v)$ and $Y = \overline{N}(v) = V(H)$. Then
\[
\begin{array}{ll}
\tau(G) & = c(G) + \overline{c}(G) - n(G)\\
& \geq (1 + c(H)) + (\overline{c}(H) + |X|) - (|X| + |Y| + 1)\\
& = c(H) + \overline{c}(H) - |Y| = \tau(H).
\end{array}
\]
\end{proof}

\begin{lemma}
\label{L33}
If $\tau(G) \leq 2$ then $\omega(G) + \alpha(G) \geq n(G)$.
\end{lemma}

\begin{proof}
We use induction on $n(G)$. Let $G$ be a graph, $n(G)=n$, $\tau(G) \leq 2$, so $c(G) + \overline{c}(G) \leq n+2$. By Theorem \ref{T28}, $d(G) + \overline{d}(G) \leq 3$, so $d(G) = 1$ or $\overline{d}(G)=1$. We may assume that $d(G)=1$, i.e., some vertex $v$ is in a unique maximal clique; i.e., $N(v)$ is a clique. We may assume that $N(v) \neq \varnothing \neq \overline{N}(v)$. Let $Z = N(v)$, $W = \overline{N}(v)$; then $n = |W| + |Z| + 1$. Let $H = G[W]$; Then $\tau(H) \leq 2$ by Lemma \ref{L32}. By the inductive hypothesis, $\omega(H) + \alpha(H) \geq |W|$. By Lemmas \ref{L31} and \ref{L29}, $H$ is $K_1$-like or $P_4$-like or $C_4$-like or $\overline{C}_4$-like. Since $A = Z \cup \{v\}$ is a clique which meets every maximal anticlique, and since $\overline{D}(v) = \overline{c}(H)$, we have
\begin{equation}
\label{L33 eqn a}
\overline{c}(G) = \sum_{a \in A}\overline{D}(a) = \overline{c}(H) + |Z| + \sum_{z \in Z}(\overline{D}(z) - 1).
\end{equation}

\item[\emph{\textbf{Case 1.} $H$ is $K_1$-like.}]

So $W = X \cup Y$ where $X$ is a clique, $Y$ is an anticlique, and $X \cap Y = \{p\}$. Then $n = |X| + |Y| + |Z|$ and $\overline{c}(H) = |X|$, so 
\begin{equation}
\label{L33 eqn b}
\overline{c}(G) = |X| + |Z| + \sum_{z \in Z}(\overline{D}(z) - 1) = n - |Y| + \sum_{z \in Z}(\overline{D}(z) - 1).
\end{equation} 

Let $B = Y \cup \{v\}$ and let $m$ be the number of maximal cliques disjoint from $B$. Since $B$ is an anticlique,
\begin{equation}
\label{L33 eqn c}
c(G) = m + \sum_{b \in B}D(b) \geq m + 1 + |Y| + (D(p)-1).
\end{equation}

Adding \ref{L33 eqn b} and \ref{L33 eqn c} we get
$$c(G) + \overline{c}(G) \geq n + 1 + m + (D(p)-1) + \sum_{z \in Z}(\overline{D}(z)-1).$$

Since $c(G) + \overline{c}(G) \leq n+2$, it follows that
\begin{equation}
\label{L33 eqn d}
m + (D(p)-1) + \sum_{z \in Z}(\overline{D}(z)-1) \leq 1.
\end{equation}

Since $\alpha(G) = \alpha(H) + 1 = |Y| +1$, we have to show that $\omega(G) \geq |X| + |Z| - 1$, i.e., that $\omega(G) \geq |Z| + |X_0|$ where $X_0 = X \setminus \{p\}$. If $Z \cup X_0$ is a clique we're done, so we may assume that some vertex $z_0 \in Z$ is not joined to some vertex $x_0 \in X_0$. Then either $D(p)\geq 2$ (if $z_0$ is joined to $p$) or $\overline{D}(z_0) \geq 2$ (otherwise), so $m=0$, i.e., every maximal clique meets $B$. Hence it will suffice to show that $Z \cup (X_0 \setminus \{x_0\})$ is a clique, since, being disjoint from $B$, it's not a maximal clique.

\item[\emph{\textbf{Case 1A.} $z_0$ is joined to $p$.}]

So $D(p) = 2$, and $\overline{D}(z)=1$ for all $z \in Z$. Then $z_0$ is joined to every vertex in $X_0 \setminus \{x_0\}$, since $z_0$ is not joined to $x_0$, and $\overline{D}(z_0) = 1$. Likewise, if $z \in Z$ and $z$ is not joined to $p$, then $z$ is joined to every vertex in $X_0$, since $\overline{D}(z) = 1$. Finally, suppose there are vertices $z_1 \in Z \setminus \{z_0\}$ and $x_1 \in X_0 \setminus \{x_0\}$ such that $z_1$ is joined to $p$ and not to $x_1$. Then the cliques $\{p, z_0, z_1\}$, $\{p, z_0, x_1\}$, and $X$ extend to three different maximal cliques, contradicting $D(p) = 2$.

\item[\emph{\textbf{Case 1B.} $z_0$ is not joined to $p$.}]

So $\overline{D}(z_0) = 2$, $\overline{D}(z)=1$ for all $z \in Z \setminus \{z_0\}$, and $D(p)=1$. Suppose $z \in Z \setminus \{z_0\}$; if $z$ is joined to $p$ then $z$ is joined to every vertex in $X_0$ because $D(p) = 1$; if $z$ is not joined to $p$, then $z$ is joined to every vertex in $X_0$ because $\overline{D}(z) = 1$. Finally, $z_0$ is joined to every vertex in $X_0 \setminus \{x_0\}$ because otherwise we would have $\overline{D}(z_0) \geq 3$. 

\vspace{10px}

In Cases 2 -- 4 we have $\omega(H) + \alpha(H) = n(H) = |W|$ by Lemma \ref{L30}, and we have to show that $\omega(G) + \alpha(G) \geq n$. Since $n = |Z| + |W| + 1 = |Z| + \omega(H) + \alpha(H) + 1 = |Z| + \omega(H) + \alpha(G)$, we have to show that $\omega(G) \geq |Z| + \omega(H)$.

\item[\emph{\textbf{Case 2.} H is $P_4$-like.}]

So $W = X \cup Y$ where $X$ is a maximal clique in $H$, $Y$ is a maximal anticlique in $H$, $X \cap Y = \varnothing$, and $n = |X| + |Y| + |Z| + 1$. Then $\overline{c}(H) = |X| + 1$, so (\ref{L33 eqn a}) becomes
\begin{equation}
\label{L33 eqn e}
\overline{c}(G) = |X| + 1 + |Z| + \sum_{z \in Z}(\overline{D}(z)-1).
\end{equation}

The set $B = Y \cup \{v\}$ is an anticlique. Let $m$ be the number of maximal cliques disjoint from $B$; $m \geq 1$, since $X$ is contained in a maximal clique which is disjoint from $B$. Then we have 
\begin{equation}
\label{L33 eqn f}
c(G) = m + \sum_{b \in B}D(b) \geq 2 + |Y| + \sum_{y \in Y}(D(y)-1).
\end{equation}

Adding (\ref{L33 eqn e}) and (\ref{L33 eqn f}) and recalling that $n = |X| + |Y| + |Z| + 1$, we get

$$c(G) + \overline{c}(G) \geq n+2 + \sum_{y \in Y}(D(y)-1) + \sum_{z \in Z}(\overline{D}(z) - 1).$$

Since $c(G) + \overline{c}(G) \leq n+2$, it follows that $D(y) = 1$ for all $y \in Y$, and $\overline{D}(z)=1$ for all $z \in Z$. We have to show that $\omega(G) \geq |Z| + \omega(H)$, i.e., $\omega(G) \geq |Z| + |X|$. It will suffice to show that $Z \cup X$ is a clique. Assume for a contradiction that some vertex $z_0 \in Z$ is not joined to some vertex $x_0 \in X$. Since $Y$ is a maximal anticlique in $H$, $x_0$ is joined to some vertex $y_0 \in Y$. But then either $D(y_0) \geq 2$ (if $y_0$ is joined to $z_0$) or else $\overline{D}(z_0) \geq 2$ (otherwise); either way we have a contradiction.

\item[\emph{\textbf{Case 3.} H is $C_4$-like.}]

So there is a $4$-element set $U = \{p, q, r, s\} \subseteq W$ such that $pq, qr, rs, ps \in E(G)$ while $pr, qs \notin E(G)$; and $W \setminus U = X \cup Y$ where $X$ is a clique and $Y$ is an anticlique, each vertex in $X$ is joined to $p$ and $q$ and at least one more vertex in $U$, and no vertex in $Y$ is joined to any vertex in $U$. Then $n = |X| + |Y| + |Z| + 5$, $\overline{c}(H) = \omega(H) = |X| + 2$, and (\ref{L33 eqn a}) becomes
\begin{equation}
\label{L33 eqn g}
\overline{c}(G) = |X| + 2 + |Z| + \sum_{z \in Z}(\overline{D}(z) - 1).
\end{equation}

Since $B = Y \cup \{p, r, v\}$ is an anticlique and $D(r) \geq 2$,
\begin{equation}
\label{L33 eqn h}
c(G) \geq \sum_{b \in B}D(b) \geq |Y| + D(p) + 3.
\end{equation}

Adding (\ref{L33 eqn g}) and (\ref{L33 eqn h}), with $n = |X| + |Y| + |Z| + 5$, we get
$$c(G) + \overline{c}(G) \geq n + D(p) + \sum_{z \in Z}(\overline{D}(z)-1).$$

Since $c(G) + \overline{c}(G) \leq n+2$, we must have $D(p) = 2$ and $\overline{D}(z)=1$ for all $z \in Z$. A similar argument, using the anticlique $B' = Y \cup \{q, s, v\}$ instead of $B$, shows that $D(q) = 2$. We have to show that $\omega(G) \geq |Z| + \omega(H)$, i.e., $\omega(G) \geq |Z| + |X| + 2$.

\item[\emph{Claim 1. $Z \cup X$ is a clique.}]
\item \begin{proof}
Suppose $z \in Z$, $x \in X$, $z$ not joined to $x$. Then $z$ must be joined to $p$ and $q$, since $\overline{D}(z) = 1$. But then the cliques $\{p, q, z\}$, $\{p, q, x\}$, $\{p, s\}$ extend to three different maximal cliques, contradicting $D(p) = 2$.
\end{proof}

\item[\emph{Claim 2. Either $Z \cup X \cup \{q\}$ or $Z \cup X \cup \{s\}$ is a clique.}]
\item \begin{proof}
If neither $Z \cup X \cup \{q\}$ nor $Z \cup X \cup \{s\}$ is a clique, then there are vertices $t_1, t_2 \in Z \cup X$ (not necessarily distinct) such that $t_1$ is not joined to $q$ and $t_2$ is not joined to $s$. Then $t_1$ and $t_2$ are joined to $p$, because every vertex in $X$ is joined to $p$, and if $t_i \in Z$ then $t_i$ must be joined to $p$ because $\overline{D}(t_i) = 1$. But then $\{p, q\}$, $\{p, s\}$, and $\{p, t_1, t_2\}$ are cliques extending to three different maximal cliques, contradicting $D(p) = 2$.
\end{proof}

\item[\emph{Claim 3. Either $Z \cup X \cup \{p\}$ or $Z \cup X \cup \{r\}$ is a clique.}]
\item \begin{proof}
Similar to Claim 2, using $D(q)=2$.
\end{proof}

It follows from Claims 2 and 3 that there are two adjacent vertices $u, w \in U$ such that $Z \cup X \cup \{u, w\}$ is a clique, so $\omega(G) \geq |Z| + |X| + 2$.

\item[\emph{\textbf{Case 4.} $H$ is $\overline{C}_4$-like.}]
So there is a $4$-element set $U = \{p, q, r, s\} \subseteq W$ such that $pr, qs \in E(G)$ while $pq, qr, rs, ps \notin E(G)$;  and $W \setminus U = X \cup Y$ where $X$ is an anticlique and $Y$ is a clique, no vertex of $X$ is joined to $p$ or $q$ or to more than one vertex in $U$, and every vertex in $Y$ is joined to every vertex in $U$. Then $n = |X| + |Y| + |Z| + 5$, and $\overline{c}(H) = |Y| + 4$ so (\ref{L33 eqn a}) becomes
\begin{equation}
\label{L33 eqn i}
\overline{c}(G) = |Y| + 4 + |Z| + \sum_{z \in Z}(\overline{D}(z) - 1).
\end{equation}

We have to show that $\omega(G) \geq |Z| + \omega(H)$, i.e., that $\omega(G) \geq |Z| + |Y| + 2$. The set $B = X \cup \{p, q, v\}$ is an anticlique. Let $m$ be the number of maximal cliques disjoint from $B$. Then

\begin{equation}
\label{L33 eqn j}
c(G) = m + \sum_{b \in B}D(b) \geq m + |X| + D(p) + 2.
\end{equation}

Adding (\ref{L33 eqn i}) and (\ref{L33 eqn j}), since $n = |X| + |Y| + |Z| + 5$, we get

$$c(G) + \overline{c}(G) \geq n + 1 + D(p) + m + \sum_{z \in Z}(\overline{D}(z)-1).$$

But $c(G) + \overline{c}(G) \leq n+2$, so $D(p) = 1$, and $\overline{D}(z) = 1$ for all $z \in Z$, and $m = 0$, i.e., every maximal clique meets $B$.

Now $Z \cup Y \cup \{r\}$ is a clique, since $Z$ and $Y \cup \{r\}$ are cliques, and if some vertex $z \in Z$ were not joined to some vertex $w \in Y \cup \{r\}$, then, since $p$ is joined to $w$, we would have either $D(p) \geq 2$ (if $z$ is joined to $p$) or else $\overline{D}(z) \geq 2$ (otherwise). So $Z \cup Y \cup \{r\}$ is a clique of size $|Z| + |Y| + 1$, but it can't be a maximal clique since it's disjoint from $B$, so $\omega(G) \geq |Z| + |Y| + 2$.
\end{proof}

\begin{theorem}
\label{T34}
Let $G$ be a graph of order $n$. 
  \[
\begin{array}{ll}
 (a) & c(G) + \overline{c}(G) = n+1$ if and only if $G$ is $K_1$-like.$\\
 (b) & c(G) + \overline{c}(G) = n+2$ if and only if $G$ is $P_4$-like, $C_4$-like or $\overline{C}_4$-like.$\\
\end{array}
  \]
\end{theorem}

\begin{proof}
By Lemmas \ref{L29} and \ref{L30}, $c(G) + \overline{c}(G) = n+1$ if $G$ is $K_1$-like, while $c(G) + \overline{c}(G) = n+2$ if $G$ is $P_4$-like or $C_4$-like or $\overline{C}_4$-like. On the other hand, if $c(G) + \overline{c}(G) \leq n+2$, then by Lemmas \ref{L33} and \ref{L31} and \ref{L29}, $G$ is $K_1$-like or $P_4$-like or $C_4$-like or $\overline{C}_4$-like.
\end{proof}

Only part (b) of Theorem \ref{T34} is new; part (a) was proved in \cite{GK71}.

\begin{corollary}
\label{C35}
Let $G$ be a graph of order $n$.
  \[
\begin{array}{ll}
 (a) & c(G) + \overline{c}(G) \geq n + 1$.$\\
 (b) & $If $ \omega(G) + \alpha(G) \leq n$ then $c(G) + \overline{c}(G) \geq n+2$.$\\
 (c) & $If $ \omega(G) + \alpha(G) < n$ then $c(G) + \overline{c}(G) \geq n+3$.$\\
\end{array}
  \]
\end{corollary}

\begin{proof}
This follows from Theorem \ref{T28}, Theorem \ref{T34}(a), and Lemma \ref{L33}.
\end{proof}

Only part (c) of Corollary \ref{C35} is new; parts (a) and (b) were proved in \cite{GK71}.

\section*{Acknowledgement}
We would like to thank Ronald L. Graham and Jack van Lint for permission to include their Theorem \ref{T18} and the corresponding portion of Theorem \ref{T21}. Theorem \ref{T18} is also used in the proof of Theorem \ref{T27}.

\end{document}